\documentclass[11pt,reqno]{amsart}

\theoremstyle{plain}

\theoremstyle{definition}

\def\Q{\mathbb Q}   
\def\N{\mathbb N}   
\def\R{\mathbb R}    
\def\C{\mathbb C}
\def\F{\mathbb F}

\def\Z{\mathbb Z}
\def\P{\mathbb P}

\def\m{{\mathfrak m}}
\def\rank{\mbox{rank}}
\def\ker{\mbox{ker}}
\def\coker{\mbox{coker}}
\def\L{{\mathcal H}}
\def\M{{\mathcal M}}
\def\lL{\gamma}

\begin{document}
\begin{center}
{\bf HYPERSURFACE COMPLEMENTS, ALEXANDER MODULES} \\
{\bf  AND  MONODROMY}
\end{center}

\vspace{3mm}

\centerline {\bf  by  Alexandru Dimca  and Andr\'as N\'emethi }

\vskip1truecm

\section{Introduction}

\bigskip
Let $X \subset \C^{n+1}$ (resp. $V \subset \P^{n+1}$) be an algebraic hypersurface and set $M_X=\C^{n+1} \setminus X$ (resp. $M_V=\P^{n+1} \setminus V$) where we suppose $n>0$. The study of the topology of $X$, $V$ and of their complements $M_X$, $M_V$ is a classical subject going back to Zariski. In a sequence of papers Libgober has introduced and studied the Alexander invariants associated to $X$, $V$, see for instance [Li0-3].

In the affine case, let $f=0$ be a reduced equation for $X$. One can use the results on the topology of polynomial functions, see for instance [B], [ACD], [NZ], [SiT],
to study the topology of the complements $M_X$, as in the recent paper by Libgober and Tib\u ar [LiT].

By taking generic linear sections and using the (affine) Lefschetz theory,
 see for instance Hamm [H] (and 
 [Li2], [LiT], [D1] for different applications), 
one can restrict this study to hypersurfaces $X$ having only isolated singularities including at infinity, see [Li2], or in the polynomial framework, to polynomials having only isolated singularities including at infinity with respect to a compactification as in [SiT]. Simple examples show that neither of these two restricted situations is a special case of the other, hence both points of view have their advantages. However, the polynomial point of view embraces larger classes of examples due to the fact that the best compactification of a polynomial function is not usually obtained by passing from the affine space $\C^{n+1}$ to the projective space $\P^{n+1}$. This is amply explained in [LiT].

In the present paper we consider an arbitrary polynomial map $f$
(whose generic fiber is connected) and we study the Alexander invariants
of $M_X$ for any fiber $X$  of $f$. 

The article has two major messages. First, the most important qualitative 
properties of the Alexander modules (cf. \ref{3.6}, \ref{n2}, \ref{n7}
 and \ref{gen})
are completely independent of the behaviour of $f$ at infinity, or about
the special fibers. (On the other hand, for particular families of
polynomial maps with some additional information about the special fibers
or about the behaviour at infinity, one can  obtain nice vanishing or 
connectivity results; see e.g. our case of h-good polynomials below). 

The second message is that all 
the Alexander invariants of all the fibers of the polynomial $f$ 
 are closely related to the  monodromy representation of $f$. 
In fact, all the torsion parts of the Alexander modules (associated with 
all the possible fibers) can be obtained
by factorization of a unique universal Alexander module, which
is  constructed from the  monodromy representation. 
This explains nicely and conceptually all the divisibility 
properties  that have  appeared recently  in the literature connecting
 the Alexander polynomials of $M_X$
 and the characteristic polynomials of some special monodromy operators,
see [Li2] and [LiT]. 
[Note that the monodromy considered by Libgober in [Li2], section
2, is associated to a Lefschetz pencil and hence quite different
from our monodromy associated to an arbitrary polynomial.]

Nevertheless, in order to exemplify our general theory,  and also to generalize
some connectivity results already present in the literature,
we introduce the family of h-good polynomials. The family includes e.g. 
all the ``good'' polynomials considered by Neumann and Rudolph [NR], and
the polynomials with isolated singularities on the affine space and at 
infinity in the sense of Siersma-Tib\u ar [SiT]. This family of h-good
polynomials fits perfectly to the study of  Alexander invariants, and it
is our major source of examples. For different vanishing and connectivity 
results, see \ref{p23}, \ref{con}, \ref{con2} and \ref{3.6}(v).

The content of our paper is the following. 
In section 2 we establish some properties of the corresponding fundamental 
groups which basically will guide all the covering properties considered later.
Moreover, here we introduce and start to discuss the  h-good polynomials. 
In section 3 we discuss some general facts on the homology groups 
$H_*(M_X, \Z)$ concentrating on non-vanishing results for $H_n(M_X,\Z)$ 
and on $\Z$-torsion problems. 
This latter aspect was somewhat neglected recently in spite of the 
pioneering work by Libgober [Li0] and 
a famous conjecture on hyperplane arrangement Milnor fibers (see \ref{2.12}). 

In section 4 we collect some facts on (torsion) Alexander modules and prove 
one of the main results, Theorem \ref{3.6}.
In order to emphasize the parallelism of h-good polynomials with 
the case of hypersurfaces with only isolated singularities
including at infinity considered by Libgober, in some of our applications 
we recall  Libgober's results [Li2] as well. 

In the fifth  section we explain the relationship between individual monodromy
operators and Alexander modules. The two main examples, i.e. 
the monodromy at infinity and the monodromy around the fiber $X$
are discussed with special care.
These two monodromy operators have been intensively studied recently using various techniques (mixed Hodge structures, D-modules), see the references given in section 5. Via our results, all this information on the monodromy operators yields valuable information on Alexander invariants of $M_X$. Remark 
\ref{4.6} relates the homology of the cyclic coverings $M_{X,d}$ to the $d$-suspension of the polynomial $f$ and in this way the results on the Thom-Sebastiani construction in [DN2] become applicable.

Section 6 introduces into the picture not only individual monodromy operators but also the whole monodromy representation of $f$. We define two new Alexander modules associated to $f$, namely the global Alexander module $M(f)$ which can be regarded as a commutative version of the  monodromy representation,
 and a local Alexander module $M(f,b)$ associated to any fiber $X=f^{-1}(b)$.
 This module $M(f,b)$ gives a very good approximation of the classical Alexander module $H_*(M_X^c, \Z)$ of $X$ (see below for the necessary definitions). 

As a convincing example of the power of this new approach, we compute at the end the various Alexander modules for a polynomial $\C^4 \to \C$ for which a partial information on the monodromy representation is known. This examples shows in particular that the isomorphism $M(f,b) =H_n(M_X^c, \Z)$ does not always hold. 

\vspace{2mm}

We thank D. Arapura, A. Libgober, C. Sabbah
 and A. Suciu for useful discussions.

\section{Topological preliminaries, connectivity properties}

\subsection{}\label{00} 
Let $f: \C^{n+1} \to \C$ be a polynomial function with $n\geq 1$. 
It is well known that there
 is a (minimal) finite bifurcation set $B_f$ in $\C$ such that $f$ is a $C^{\infty}$-locally trivial fibration over $\C \setminus B_f$.
If $b_0 \in \C$ is not in $B_f$, then $F=f^{-1}(b_0)$ is 
called the generic fiber of $f$; otherwise 
$F_b:=f^{-1}(b)$ is  called a special fiber. 

For any $b\in \C$ we fix a sufficiently small closed  disc $D_b$ containing
$b$, and a point $b'\in \partial D_b$. We set 
$T_b:=f^{-1}(D_b)$, \ $T^*_b:=T_b\setminus f^{-1}(b)$. Sometimes, it is 
convenient to identify $f^{-1}(b') $ with the generic fiber $F$.
Then, we have the obvious inclusions $F\subset T^*_b\subset T_b$.

By a well-known deformation retract argument (see e.g. (2.3) in [DN1]),
the pair $(\C^{n+1},F)$ has the homotopy 
type of the  space $(Y,F)$ obtained by gluing all the pairs $(T_b,F)$ 
($b\in B_f$) along $F$. We denote this fact by
\begin{equation*}
(\C^{n+1},F)\sim \bigvee_{F}(T_b,F)\ \ (b\in B_f).
\tag{1}
\end{equation*}

\subsection{Proposition}\label{t1} {\em 
Let $f:\C^{n+1}\to \C$ be a polynomial map.
Then the generic fiber $F$ is  connected  if and only if 
 $\pi_1(T_b,F)$ is trivial for any $b\in B_f$.}

\begin{proof}
If $\pi_1(T_b,F)=1$ for all $b$,  then
$\tilde{H}_0(F)=H_1(\C^{n+1},F)=\oplus_{b\in B_f}H_1(T_b,F)=0$
 by (1). 
Now, assume that $F$ is connected and fix a $b\in B_f$. 
Then we have to show that  $j:\pi_1(F)\to \pi_1(T_b)$ is onto.
Since $T_b$ is smooth and $f^{-1}(b)\subset T_b$ has real codimension
two, one obtains that $i:\pi_1(T^*_b)\to \pi_1(T_b)$ is onto. 
Since $f$ restricted to $T^*_b$ is a fiber bundle,
the kernel of $f_*:\pi_1(T^*_b)\to \Z$ is $\pi_1(F)$.
Assume that $f^{-1}(b)$ has $r$ irreducible components, and
$f-b=\prod_{i=1}^rg_i^{m_i}$. Then one can construct easily elementary
loops in $T^*_b$ around the component $\{g_i=0\}$ representing $x_i\in \pi_1(
T^*_b)$  with  properties
$f_*(x_i)=m_i$ and $i(x_i)=1$. Set $m:=gcd_i\{m_i\}$.
Then a combination of the $x_i$'s provides an $x\in \pi_1(T^*_b)$
with $f_*(x)=m$ and $i(x)=1$. The point is that $m=1$ (otherwise
$f-b$ would be an $m$-power of a polynomial whose generic fiber is not
connected). The existence of such an $x$ and the surjectivity of $i$
implies the surjectivity of $j$. \end{proof}

In the next paragraphs we fix a $b\in B_f$,  and we write $X:=f^{-1}(b)$
and $M_X:=\C^{n+1}\setminus X$. For simplicity of the notations, we will assume
that $b=0$. 

\subsection{Corollary}\label{t2} {\em Assume that $F$ is connected. Then

\vspace{1mm}

(i) $ \pi_1(T^*_0)\to \pi_1(M_X)$ (induced by the inclusion) is onto.

\vspace{1mm}

(ii)  $\pi_1(F)\stackrel{i_X}{\longrightarrow}
 \pi_1(M_X)\stackrel{f_*}{\longrightarrow}\Z\to (1)$\ 
is an exact sequence (i.e. $im(i_X)=ker(f_*)$), where 
$i_X$ is induced by the inclusion  $F\subset M_X$, and $f_*$ by $f:M_X\to 
\C^*$.}
\begin{proof}
Similarly as in (1),  $M_X$  has the homotopy type of a space obtained by 
gluing $T^*_0$ and all the ``tubes'' $T_{\bar{b}}$ ($\bar{b}\in B_f\setminus 
\{0\}$) along $F$.
Then (i) follows from van Kampen theorem and from the 
surjectivity of $\pi_1(F)\to \pi_1(T_{\bar{b}})$
 for each $\bar{b}$ (cf. \ref{t1}). 
Part (ii) follows from (i) and the exact sequence 
$\pi_1(F) \to  \pi_1(T^*_0)\to\Z\to (1)$. \end{proof}

Let $p:\F\to M_X$ be the $\Z$-cyclic covering associated to the kernel 
of the morphism $f_*: \pi _1(M_X) \to\Z$. The notation $\F$ is chosen because 

(i) $\F$ is the
homotopy fiber of $f:M_X \to \C^*$, regarded as a homotopy fibration; and

(ii) in many cases the topology of $\F$ is a good approximation for the topology of $F$ (see e.g. the connectivity results below). 

Fix a base-point $*\in \F$ with $p(*)\in F$. Since $f(F)$ is a point, there  
is a natural section $s:F\to \F$ of $p$ above $F$ with $s(p(*))=*$.
In particular, we can regard $F$ as a subspace of $\F$. 

\subsection{Corollary}\label{t3} {\em Assume that $F$ is connected. Then 
$s_*:\pi_1(F)\to \pi_1(\F)$ is onto, or equivalently,
$\pi_1(\F,F)$ is trivial.}
\begin{proof}
Compare the exact sequences $(1)\to \pi_1(F)\to \pi_1(T^*_0)\to \Z\to (1)$
and $(1)\to \pi_1(\F)\to \pi_1(M_X)\to \Z\to (1)$ via \ref{t2}.
\end{proof}

Fix an  orientation  of $S^1$, and consider a
smooth loop $\lL:S^1\to \C\setminus B_f$. 
Denote by $q:\lL^{-1}(f)\to S^1$ the pull-back of $f$ by $\lL$,
i.e. $\lL^{-1}(f)=\{(t,x)\in S^1\times \C^{n+1}:\ \lL(t)=f(x)\}$,
and $q(t,x)=t$. 

\subsection{Corollary}\label{t4} {\em Assume that $\lL_*:\pi_1(S^1)\to
\pi_1(\C^*)$ (i.e. $\lL_*:\Z\to\Z$) is multiplication by an integer $\ell$.
Then one has the following commutative diagram with all the lines
and columns exact:}
\[
\begin{array}{ccccccc}
(1)  \to & \pi_1(F) & \to & \pi_1(\lL^{-1}(f)) & \to & \Z & \to  (1) \\
  & \downarrow &  & \downarrow & & \downarrow &  \\
(1)  \to & \pi_1(\F) & \to & \pi_1(M_X) & \to & \Z  &\to (1) \\
  & \downarrow &  & \downarrow & & \downarrow &  \\
  & (1)  & \to & \Z/\ell\Z  & \to & \Z/\ell\Z & \to  (1) \\
  & &  & \downarrow & & \downarrow & \\
  & &  & (1) & & (1)  &  
\end{array}\]
\begin{proof}
The first two lines are the homotopy exact sequences 
of the corresponding fibrations. Then use \ref{t3}. \end{proof}

Sometimes it is convenient to work with special polynomials with nice 
behaviour  around the special fibers or at infinity.
First, we recall the definition of Neumann and Rudolph of good polynomials
[NR].   A fiber $f^{-1}(b)$ is called ``regular at infinity''
if there exist
a small disc $D$ containing $b$ and a compact set $K$ such that 
the restriction of $f$ to $f^{-1}(D)\setminus K$ is a  trivial 
$C^{\infty}$-fibration. The polynomial $f$ is called 
good (or ``topologically good'') if all its fibers are regular at infinity. 

For example,
the tame polynomials introduced by Broughton [B], 
the larger class of M-tame polynomials introduced by N\'emethi-Zaharia [NZ] 
are good. We recall that any fiber of a good polynomial is a 
bouquet of spheres $S^n$, $B_f$ is the set of critical values of $f$,
for any $b\in B_f$ the ``tube'' $T_b$ has the 
homotopy type of $f^{-1}(b)$, and $f^{-1}(b)$
(homotopically)  is obtained from 
$F$ by attaching some cells of dimension $n+1$. 

For the purpose of the present paper it is enough (and it is more natural)
a much weaker assumption.

\subsection{Definition}\label{d22} 
The polynomial $f$ is called 
``homotopically  good'' (h-good)  if for any $b\in B_f$,  the pair 
$(T_b,F)$ is $n$-connected. 

\vspace{2mm}

From the above discussion it follows easily that all the good polynomials are
h-good. Another example is provided by 
the polynomials with isolated singularities
on the affine space and at infinity in the sense of Siersma-Tib\u ar [SiT]
(see p.776 [{\em loc.\,cit.}]). 

In general, for an arbitrary polynomial, it is much easier to handle
the properties of the generic fiber and the ``tubes'' $T_b$ than
the properties of the special  fibers (see e.g. [DN2]). One of the 
advantages of the above definition \ref{d22} is that it requires
information only  about $F$ and $T_b$'s. (Conversely, this fact also
explains that for h-good polynomials one can say very little about the 
special fibers. E.g. the special fibers of  h-good polynomials,
in general,  are not even  reduced, as it happens  e.g. 
for $f(x,y)=x^2y$.
For a non-trivial example of a  h-good polynomial which has non-isolated 
singularities, see the polynomial $f_{d,a}$
constructed by tom Dieck and Petrie, cf. [D1], p.175.)

The second advantage of  definition \ref{d22} is that, in fact,  it is almost
homological. Indeed, for $n=1$, $f$ is h-good iff 
$F$ is connected (by \ref{t1}); for $n>1$ the polynomial $f$ is  h-good
iff $F$ is simply-connected and $H_q(T_b,F,\Z)=0$ for all $b$ and $q\leq n$.
This second statement follows from \ref{t1}, the next proposition \ref{p23},
 and the relative Hurewicz isomorphism theorem (see e.g. [S], p.397). 

\vspace{2mm}

The above examples and comment show that 
we cannot expect that the h-good polynomials
 will share all the properties of the ``good''  ones. However, 
the next result will recover one of the most important properties.

\subsection{Proposition.}\label{p23} {\em Assume that $f$ is a h-good 
polynomial. Then its generic fiber $F$ has the homotopy type of a bouquet
of spheres $S^n$.}

\begin{proof} 
By \ref{00}(1), 
$H_q(\C^{n+1},F)=0$ for $q\leq n$, hence $\tilde{H}_q(F)=0$ for 
$q\leq n-1$. This already proves the statement for $n=1$. 
Next, we have to show that $\pi_1(F)=(1)$, provided that $n\geq 2$. 
The connectivity assumption assures that for each $b\in B_f$, 
$\pi_1(F)\to \pi_1(T_b)$ is an isomorphism. We denote all these 
fundamental groups by $G$. If the cardinality $|B_f|$ of $B_f$ is one,
then this implies that $\pi_1(F)=G$ is trivial
 (since in this case,  $T_b\sim \C^{n+1}$). If $|B_f|>1$,
then by van Kampen theorem, applied for
$Y=\vee_F(T_b)$ (cf. (1)), 
and induction over $|B_f|$, we get $\pi_1(Y)=G$. 
But again by \ref{00}(1), $\pi_1(Y)=(1)$.
Since $F$ has the homotopy type of a finite
$n$-dimensional 
 CW complex, the result follows by Whitehead theorem. \end{proof}

\subsection{Remark}\label{2.4} (\ref{p23})  can be compared with the 
following classical result of L\^e [L\^e]. 

\noindent  {\it For any projective hypersurface $V$ and a generic
hyperplane $H$, the affine hypersurface $X=V \setminus H$ is homotopy 
equivalent to a  bouquet of spheres $S^{n}$.}
[For the computation of the number of spheres in this bouquet, in terms of the 
degree of gradient mappings,  see [DPp].]

\subsection{}\label{fin}
Finally, we compare $F$ and $\F$. Since $\F$ is a cyclic covering of $M_X$,
and $M_X$ has the homotopy type of a finite CW complex of dimension $\leq 
(n+1)$, one has the general fact:
\begin{equation*}
\mbox{$H_m(\F)=0$ for $m>n+1$ and $H_{n+1}(\F,\Z)$ has no $\Z$-torsion.}
\tag{2}
\end{equation*}
But if $f$ is h-good, we can say more (cf. also with \ref{t3}).
With a choice of the base points, we again embed $F$ into $\F$ via the
section  $s$.

\subsection{Proposition}\label{con} {\em 
If $f$ is h-good then the pair $(\F,F)$ is $n$-connected.
Therefore, $\F$ is $(n-1)$-connected, and $s_n:H_n(F,\Z)\to
H_n(\F,\Z)$ is onto. In particular, for $n>1$, by Hurewicz theorem
and the homotopy exact sequence, one has 
$H_n(\F,\Z)=\pi_n(\F)=\pi_n(M_X)$.}
\begin{proof}
Notice (cf. \ref{00}(1) and the proof of \ref{t2}) that $M_X$ has the homotopy
type of a space obtained by gluing to $T^*_0$ along $F$ all the tubes 
$T_{\bar{b}}$ for $\bar{b}\in B_f\setminus \{0\}$. Moreover, $p^{-1}(T^*_0)$
has the homotopy type of $F$, and its embedding into $\F$ is homotopically
equivalent to the embedding $s:F\to \F$. Therefore, by excision:
$$H_q(\F,F)=H_q(p^{-1}(T^*_0\vee_F(T_{\bar{b}})),p^{-1}(T^*_0))=
\oplus_{\bar{b}}H_q(p^{-1}(T_{\bar{b}}),p^{-1}(F)),$$
where $\bar{b}$ runs over $B_f\setminus\{0\}$. 
But $(p^{-1}(T_{\bar{b}}),p^{-1}(F))=\Z\times (T_{\bar{b}},F)$, hence
$H_q(\F,F)=0$ for $q\leq n$. Hence, the connectivity follows from this,
\ref{t3}, \ref{p23},  and the relative Hurewicz isomorphism theorem. 
Finally, the connectivity of $\F$ follows from the connectivity of $F$, 
cf. \ref{p23}. \end{proof}
\noindent 
The same proof (but neglecting $p$), and the Wang exact sequence
of $T^*_0,$ gives:

\subsection{Corollary}\label{con2} {\em If $f$ is h-good, then
the pair $(M_X,T^*_0)$ is $n$-connected. In  particular,
$H_q(M_X,\Z)=0$ for $1<q<n$.  (In fact, by \ref{con}, the cyclic covering
of $M_X$ is $(n-1)$-connected, hence $\pi_q(M_X)=0$ for 
$1<q<n$ as well.)}

For special cases of this connectivity results, see also  [Li2] and [LiT].

\section{Preliminaries about $H_*(M_X,\Z)$}

\noindent Let  $X$ be a hypersurface in $\C^{n+1}$ and $M_X$ be 
its complement.
The goal of this section is to list some properties of the integral
homology of $M_X$, with an extra emphasis on the torsion part
and the ``interesting part'' $H_n(M_X,\Z)$. 
Moreover, we present some constructions which generate examples 
with non-trivial ``interesting part''. 
Additionally, sometimes we compare the properties of $H_*(M_X)$ with 
homological properties of hypersurfaces $X$. 

 We start with the case when $X$ is a (generic or special)
fiber of a polynomial $f$.

\subsection{Fact}\label{2.1}[LiT] 
 {\it If $F$ is the generic fiber of an arbitrary
 polynomial $f$, then $M_F$ 
has the homotopy type of a join $S^1 \bigvee S(F)$, 
where $S(F)$ denotes the suspension of $F$.  In particular,}
\begin{equation*}
 {\tilde H}_k(M_F)={\tilde H}_k(S^1)\oplus   {\tilde H}_{k-1}(F)
\ \ \ \mbox{{\em for any $k$}}.
\tag{1}
\end{equation*}

\vspace{2mm}

In fact,  the result \ref{2.1}(1)
 holds for any smooth $X$, as follows from the associated Gysin sequence, 
see [D1], p.46.

The similar result (i.e. the analog of \ref{2.1}(1))
  for homotopy groups is definitely false; consider 
for example $\pi_1$, or (for instance) 
[S], p.419, exercise B6, for a reason.

\subsection{Example}\label{2.3} \ref{2.1}(1)
 is false for special fibers, even for very simple polynomials. 
Let $f=x_0^2+...+x_n^2$ and $X=f^{-1}(0)$. The Wang sequence
 of the global Milnor fibration
$F \to M_X \to \C^*$ (see [D1], p.71-74) and 
the fact that the corresponding monodromy operator $T$ acting on 
$H_n(F,\Z)=\Z$ is $(-1)^{n+1}Id$,  implies the following:

(i) for $n=2m+1$ odd, $H^*(M_X)=H^*(S^1 \times S^n)$. In fact,
the monodromy is isotopic to the identity and hence we have a diffeomorphism $M_X = \C^* \times F$. This implies
that $M_X$ has the homotopy type of $S^1 \times S^n$ as claimed in [Li2], 
Remark (1.3);

(ii) for $n=2m>0$ even, $H_n(M_X,\Z)=\Z/2\Z$. In particular,
$M_X$ is not homotopy equivalent to the product $S^1 \times S^n$ 
(contrary to the claim in  [Li2], Remark (1.3)). 

\vspace{2mm}

(\ref{2.1}) has the following consequence: when $X$ is the generic fiber 
of a h-good polynomial then $H_m(M_X,\Z)=0$ for $1< m \leq n$, and in fact 
$H_*(M_X,\Z)$ is torsion free (cf. \ref{p23}).  (This can be compared with 
L\^e's result \ref{2.4}, which  shows that 
{\em generically} an  affine hypersurface $X$ has no torsion in homology.)

More generally, it was shown in [Li2] that 
 when $X$ has isolated singularities 
including at infinity,  then  $H_m(M_X)=0$ for $1<m<n$.
The same statement holds for the special fiber $X$  of a 
h-good polynomial by \ref{con2} (cf. also with [LiT]). 
Hence, in both cases, the first interesting 
homology group occurs in degree $n$.

\smallskip

We describe now three constructions which provide in a 
  systematic way  hypersurfaces $X$ with $H_n(M_X, \Z) \not= 0$. 

Below  $V$ denotes the projective closure of $X$ and $H$ the hyperplane at 
infinity.

\subsection{The first construction (using duality)}\label{2.7}
Following  [Li2], section 1, we consider the isomorphism
\begin{equation*}
 H_m(M_X)=H^{2n-m+1}(V,V \cap H)\ \ \ \mbox{for all $m$.}
\tag{2}
\end{equation*}
 Assume that $V$ and  $V \cap H$ have only isolated singularities
(this is exactly the condition on $X$ to have isolated singularities including at infinity in [Li2]).
The exact sequence
\begin{equation*}
H^n(V) \to H^n(V \cap H) \to H^{n+1}(V,V \cap H) \to H^{n+1}(V) 
\to H^{n+1}(V \cap H)
\tag{3}
\end{equation*}
and the isomorphism  $H^{n+1}(V \cap H)=H^{n+1}(\P^{n-1})$ (cf. [D1], p.161)
imply the following.

\smallskip

\noindent 
{\it Assume that $V$ and \,  $V \cap H$ \, have only isolated singularities and 
that $H^{n+1}(V,\Z) \not=H^{n+1}(\P^{n-1},\Z)$. Then $H_n(M_X, \Z) \not= 0$. 
In particular, the corresponding affine hypersurface is not the generic fiber of a h-good polynomial.}

\subsection{Example}\label{2.8}  Let $V$ be a cubic surface in $\P^3$ having two singularities, one of type $A_1$ and the other of type $A_5$. First we take $H$ a generic plane.
In this case, using [D1], p.165
we see that the exact sequence \ref{2.7}(3)  becomes
$$\Z \to \Z \to H_2(M_X,\Z) \to \Z/2\Z \to 0$$
where the first morphism is multiplication by $deg(V)=3$. It follows that
$H_2(M_X,\Z)=\Z/6\Z$, $X$ is singular and the associated polynomial $f$ is tame.

Secondly, take $H$ to be any plane passing through the 2 singularities on $V$.
Then the associated $X$ is smooth, but by \ref{2.7},
 $X$ is not the generic fiber of a good polynomial.

\subsection{The second construction (using finite
cyclic coverings and defect)}\label{defect} The second approach uses 
$M_{X,e}$, the cyclic covering of $M_X$ of degree $e$ when $X=f^{-1}(0)$
(cf. also with  \ref{3.4}(II) and \ref{remarks}(4)). 
It is clear that we can take
$$M_{X,e}=\{(x,u) \in \C^{n+1} \times \C^* | f(x)-u^e=0 \}.$$
In some cases we can get a useful approximation of $M_{X,e}$ as follows.
Fix a system of positive integer weights $w=(w_0,...,w_n)$, and let
$e$ be the top degree term in $f$ with respect to $w$.
Introduce a new variable $t$ of weight 1 and let ${\tilde f}(x,t)$ be the 
homogenization of $f$ with respect to the weights $(w,1)$. Consider the 
affine Milnor fiber $F': {\tilde f}(x,t)   =1$, which is a smooth 
hypersurface in $\C^{n+2}$.
One has an embedding $j: M_{X,e} \to F'$ given by
$$j(x,u)=(u^{-1} * x, u^{-1}),$$
where $*$ denotes the multiplication associated to the system of weights $w$.
The complement $F' \setminus j(M_{X,e})$ is characterized by $\{t=0\}$,
hence it  can be identified with
 the affine Milnor fiber 
$F_e:f_e(x)=1$ (considered din $\C^{n+1}$)
defined by the top homogeneous component $f_e$ of $f$. 
If $f_e$ defines an isolated singularity at the origin, then $f_e$ is a good 
polynomial, $F_e$ is $(n-1)$-connected
and the Gysin sequence of the divisor $F_e$ implies the isomorphisms
\begin{equation*}
j_*:H_k(M_{X,e}) \to H_k(F')\ \ \mbox{for $1<k<n+1$}.
\tag{4}
\end{equation*}
Note that under this isomorphism the action of the natural generator of the  covering transformation group on $M_{X,e}$ corresponds to multiplication by $exp(-2 \pi i/e)$ (of all the coordinates) on $F'$. 
Moreover, notice that $H_k(M_X,\Q)$ is isomorphic to the group of invariants of
$H_k(M_{X,e},\Q)$ with respect to this action. 

Notice that the dimension of $H_n(F')$ is closely related to the  {\em 
defect} associated with the singular points of the projective 
hypersurface $\tilde{f}=0 $ considered
in $\P^{n+1}$ (for details, see e.g. [D1]). Hence, this construction 
emphasizes the connection  between $H_n(M_X)$ and superabundance 
properties. 

For more information on the homology of $M_{X,e}$,  see also \ref{remarks}(4)
and  \ref{4.6}  below. 

\subsection{Example}\label{2.10}  Let $f:\C^4 \to \C$ be given by 
$f(x,y,u,v)=x^3+y^3+xy-u^3-v^3-uv$. Set $X=\{f=0\}$. Then $f$ is a tame 
polynomial and $X$ has 10 nodes. 
It follows from [D1], p.208-209, (with all the weights $w_i=1$)
that $\dim H_3(F',\C)=5$ and the 
 multiplication by $exp(-2 \pi i/3)$ on $F'$ induces the trivial action on 
$H_3(F',\C)$. It follows that $\dim H_3(M_X,\C)=\dim H_3(M_{X,3},\C)=
\dim H_3(F',\C)=5$.

\subsection{The third construction (``counting'' the Milnor numbers)}\label{Mil} For simplicity we will assume that $f$ is a (topologically)
{\em good } polynomial. 
As above, $F$ is its generic fiber and $X$ a special fiber.
Assume that $X$ has $n_X$ singular points with Milnor fibers 
$F_i $ and Milnor numbers $\mu_i$ for  $1\leq i\leq n_X$. For each $i$,
let $\mu_{0,i}$ denote the  rank of $H_n(\partial F_i)$. Set $\mu_X:=
\sum_i\mu_i$ and $\mu_{0,X}:= \sum_i\mu_{0,i}$. 
Below, $\oplus_X$ means $\oplus_{i=1}^{n_X}$. 

Finally, let
$\mu$ be the sum of all the Milnor numbers of the singularities of $f$
(situated on all the singular fibers). With these notations, one has:

\vspace{2mm}

{\em 
\noindent
The groups $H_q(M_X)$ ($q=n, n+1$) are inserted in the following 
commutative diagram:

\[
\begin{array}{ccccccccccc}
0 & \to & \oplus_XH_n(\partial F_i) & 
\to & \oplus_XH_n(F_i) &
\to & \oplus_XH_n(F_i,\partial F_i) &
\to & \oplus_X\tilde{H}_{n-1}(\partial F_i) & 
\to & 0\\
 & & \downarrow &  & \downarrow & & \downarrow & & \downarrow & \\
0 & \to & H_{n+1}(M_X) &
\to & H_n(F) &
\to & \oplus_XH_n(F_i,\partial F_i) &
\to & H_{n}(M_X) &
\to & 0
\end{array}\]
where the first two vertical arrows are monomorphisms, the third is
the identity, and the last one is an epimorphism.
In particular:
$$\mu_{0,X}\geq \dim H_n(M_X) \geq \mu_X+\mu_{0,X}-\mu.$$}

Above, the first inequality is not new, and it is not very sharp
(cf. e.g. with [GN2] (2.30), 
or with [Si2] (5.4) and \S 9; it also follows from \ref{4.6}(iii) below). 

The second inequality is more interesting,
and it  can be used in two different ways. First, using the 
integers $ \mu_X$, $\mu_{0,X}$ and $\mu$, if the sum of the first two 
is strictly larger than the third, one gets a non-vanishing criteria for
$H_n(M_X)$. 
For example, in the case \ref{2.10}, $\mu=16$, $\mu_X=\mu_{0,X}=10$,
hence the second inequality reads as $\dim H_n(M_X)\geq 4$.

Similarly, if one knows $\mu_X$, $\mu_{0,X}$ and $\dim H_3(M_X)$
 about the hypersurface $X$,  then these numbers 
may  impose serious conditions about singularities
of some other singular fibers of $f$ (e.g. about their existence). 
\begin{proof}
Assume that $X=f^{-1}(b)$, and let $T_b^o$ be the interior of $T_b$. 
Clearly, $M_X$ has the homotopy
type of  $\C^{n+1}\setminus T_b^o$. Let $F$ be a fixed generic
fiber of $f$ inside of $T_b^o$. Then one can write the homological  long 
exact sequence of the pair $(\C^{n+1}\setminus F, \C^{n+1}\setminus T_b^o)$.
Notice that $H_q:=H_q(\C^{n+1}\setminus F, \C^{n+1}\setminus T_b^o)$ equals
$H_q(T_b\setminus F, \partial T_b)$ by excision. 
Let $B_i$ be a small Milnor ball
of the $i$-th singular point of $X$. Then using the ``good''-property 
of $f$  and
excision one gets $H_q=\oplus_XH_q
(B_i\setminus F_i, \partial B_i\setminus \partial F_i)$.
But this is isomorphic to $\oplus_XH_{q-1}(F_i,\partial F_i)$ by
Gysin isomorphism. Use these facts and the Gysin isomorphism
$H_n(F)\to H_{n+1}(M_F)$ to obtain the second line of the above diagram. 

Next, assume that the disc $D_b$ is sufficiently small with respect to the balls $B_i$, and  consider for any ball $B_i$ the local analog of the above
 picture; namely
the homological long exact sequence of the pair 
$(B_i\setminus F,B_i\setminus 
T_b^o)$. This homological sequence admits a natural map to the previous 
sequence induced by the inclusion. Finally, this ``local sequence''
is modified by dualities and Gysin isomorphisms. \end{proof}

\subsection{Remarks about the torsion part of $H_n(M_X,\Z)$.}\label{torsion}
In the final part of this section we discuss the relations between 
the existence/non-existence of  torsion
in the homology of $X$ and $M_X$ respectively. 
The following examples show that
these relations are not simple even for a homogeneous polynomial $f$. 

\subsection{Example}\label{2.11} Consider 
the homogeneous polynomial $f=x^2y^2+y^2z^2+z^2x^2-2xyz(x+y+z)$ defined on 
$\C^3$.
It is known that its Milnor fiber $F$ has torsion in homology, 
more precisely $H_1(F,\Z)= H_2(M_{F},\Z) = \Z/3\Z$ and 
$H_2(F,\Z)=\Z^3$ see [Li0], [DN1] and [Si].
Let $C$ be the 3-cuspidal quartic in $\P^2$ defined by $\{f=0\}$.
It satisfies $H_1(\P^2\setminus C,\Z)=\Z_4$ and 
$H_2(\P^2\setminus C,\Z)=0$ (cf. [loc. cit.]). 
Set $X=f^{-1}(0)$.  
The Gysin sequence of the fibration $\C^* \to M_{X} \to \P^2 \setminus C$ 
yields $H_1(M_X,\Z)=\Z$, $H_2(M_{X},\Z)=\Z/4\Z$ and  $H_3(M_{X})=0$.
In conclusion, both $M_F$ and $M_X$ have torsions in $H_2$, but these 
torsions are different.

\vspace{2mm}

\noindent 
If we put together the examples \ref{2.3} and \ref{2.11},
 we see that for a homogeneous polynomial $f$ and for $X=f^{-1}(0)$ 
the only case not covered is the following.

\subsection{Question}\label{2.12} 
{\em Find an example of a homogeneous polynomial $f$ such that
$M_X$ has no torsion but $F$ has torsion.}

 Even in  the case when $f$ is a product of linear forms,
 the existence of such an example is an open question. (It is known in this latter case that the hyperplane arrangement complement $M_X$ is torsion free, see [OT],
 and the corresponding Milnor fiber can be identified to a cyclic covering of
$\P^n\setminus \{f=0\}$, see [CS], [CO]).
Notice also, that if we allow $f$ to be a product of {\em powers
 of linear forms}, then A.
Suciu has examples with torsion in the homology of the associated
Milnor fiber $F$.

\smallskip

The following result gives some conditions on the possible torsions that may arise in such a case.

\subsection{Proposition}\label{2.13} {\it Assume 
that for a homogeneous polynomial $f$ of degree $d$ the complement $M_X$ has no $p$-torsion for some prime $p$ and that the $p$-torsion in a homology group $H_k(F,\Z)$ has the structure
$$\Z/p^{k_1}\Z \oplus \Z/p^{k_2}\Z \oplus ...  \oplus  \Z/p^{k_m}\Z$$
with  $m \geq 1$ and $k_1 \geq k_2 \geq ... \geq k_m \geq 1$. Then

(i) $(p-1,d)=1$ implies $m \geq 2$ and $k_1=k_2$;

(ii) $((p-1)p(p+1),d)=1$ implies $m \geq 3$ and $k_1=k_2=k_3$.  }
\begin{proof} Consider 
the Wang sequence in homology associated to the fibration
$F \to M_X \to \C^*$. The fact that $M_X$ has no $p$-torsion implies that $T-Id$ is
an isomorphism when restricted to the $p$-torsion part.

To prove (i) note that since such an isomorphism preserves the orders of the elements, unless the claim (i) above holds, we get
an induced automorphism of $\Z/p\Z$, where this latter group is regarded as the quotient of the $p$-torsion part by the subgroup of elements killed by multiplication by $p^{k_1-1}$. The same is true for the monodromy transformation $T$.

Denote by $T_p$ the induced automorphism of $\Z/p\Z$.
It follows that $T_p$ is not the identity, $T_p^d=Id$ (since $f$ is 
homogeneous of degree $d$) and $T_p^{p-1}=1$ (since $|Aut(\Z/p\Z)|=p-1$), 
a contradiction.

For (ii) use the same argument, plus the equality $|Aut((\Z/p\Z)^2)|=(p-1)^2p(p+1)$. \end{proof}

\subsection{Remark}\label{2.14} Other examples involving torsion in the homology of the special fibers of a polynomial can be obtained by suspension, see for details [DN2].

\section{Alexander Modules}

\subsection{Definitions}\label{3.1}
Let $Y$ be a connected CW-complex and $e_Y: \pi _1(Y) \to \Z$ be an epimorphism.
We denote by $Y^c$ the $\Z$-cyclic covering associated to the kernel of the morphism $e_Y$.  It follows that a generator of $\Z$ acts on $Y^c$ by a certain covering homeomorphism $h$ and all the groups $H_*(Y^c,A)$, $H^*(Y^c,A)$ and $\pi _j(Y^c) \otimes A$ for $j>1$ become in the usual way $\Lambda _A$-modules, where $\Lambda _A= A[t,t^{-1}]$, for any ring $A$. These are called the Alexander modules of the pair $(Y,e_Y)$ or simply of $Y$ when the choice of $e_Y$ is clear.

If $Z$ is a second connected CW-complex, $e_Z: \pi _1(Z) \to \Z$ an epimorphism and $\phi : Y \to Z$ is a continuous map such that the induced map at the level of $\pi _1$ is an epimorphism compatible to the two given epimorphisms $e_Y$ and $e_Z$, then $Y^c \to Y$ can be regarded as a pull-back covering obtained via $\phi$ from the covering $Z^c \to Z$.
In particular, this gives a lift $\phi ^c :Y^c \to Z^c$ which is compatible 
with the covering transformations, and hence,
 an induced morphism of $\Lambda _A$-modules, say $\phi ^c_*: H_*(Y^c,A) \to   H_*(Z^c,A)$.

If $A$ is a field then the ring $\Lambda _A$ is a PID. Hence any finite type $\Lambda _A$-module $M$ has a decomposition $M=\Lambda _A ^k \oplus (\oplus _pM_p)$, where $k$ is the rank of $M$ and the second sum is over all the prime
elements in $\Lambda _A$,  and $M_p$ denotes the $p$-torsion part of $M$.

More precisely, for each prime $p$ with $M_p\not=0$,
we have a unique decomposition

\begin{equation*}
M_p= \oplus_{i=1,\ell_p}\Lambda _A/p^{k_i}
\tag{1}
\end{equation*}
for $\ell_p>0$ and $k_1 \geq k_2 \geq .... \geq k_{\ell_p} \geq 1$.
The sequence $(k_1,k_2,...,k_{\ell_p})$ obtained in this way will be denoted by
$K(M,p)$. One can define an order relation on such sequences by saying that
$$(k_1,....,k_a) \geq (m_1,...,m_b)$$
iff $a \geq b$ and $k_i \geq m_i$ for all $1 \leq i \leq b$.

Let $\Delta _p(M)=\prod _{i=1,\ell_p}p^{k_i}$ (resp. $\Delta (M)=\prod _p \Delta _p(M)$) be the $p$-Alexander polynomial (resp. the Alexander polynomial) of the module $M$. This latter invariant $\Delta (M)$ is called the order of $M$ in [Li2]. See \ref{remarks}(1)  for a motivation of this terminology.

With this notation one has the following easy result whose proof is left to the reader.

\subsection{Lemma}\label{3.2} {\it Let $u:M \to N$ be an epimorphism of $R$-modules, where $R$ is a PID and $M$ is of finite type. Then $N$ is of finite type and for any prime $p \in R$ one has $K(M,p) \geq K(N,p)$. In particular $\Delta (N)$ divides $\Delta (M)$.}

\subsection{Example}\label{3.3}
 For $A=\C$, we will simply write $\Lambda $ instead of $
\Lambda_{\C}$. The prime elements in this case are just the linear forms 
$t-a$, for $a \in \C^*$. Moreover, a $\Lambda $-module of the form 
$H_q(Y^c,\C)$ is of finite type and torsion iff the corresponding Betti number 
$b_q(Y^c)$ is finite. If this is the case, the $(t-a)$-torsion part of 
$H_q(Y^c,\C)$ is determined (and determines) the Jordan block structure of 
the corresponding automorphism $h_q$; i.e. a Jordan block of size $m$ 
corresponding to the eigenvalue $a$ produces a summand $\Lambda /(t-a)^m$. 
The corresponding Alexander polynomial $\Delta (H_q(Y^c,\C),h_q)$ is just 
the characteristic polynomial $\Delta (h_q)$ of $h_q$.

When $A=\Q$, the corresponding prime elements are the irreducible polynomials in $\Q[t]$ different from $t$, hence they are a lot more difficult to describe.
However, the knowledge of the $\Lambda $-module structure implies easily the
$\Lambda_{\Q} $-module structure just by grouping together the polynomials $t-a$ for those $a$'s having the same minimal polynomial over $\Q$.

\subsection{The Alexander modules of $M_X$ and local systems}\label{3.4} 
Coming back to the situation (and notation) of the previous sections,
for any hypersurface  $X$ we define 
$e_X:\pi_1(M_X)\to \Z$ as follows.  In fact, we will distinguish two cases.

\vspace{2mm}

\noindent (I) \ First, assume that $X$  is an arbitrary hypersurface 
in $\C^{n+1}$ (and even if $X=f^{-1}(0)$,   we disregard  $f$). 
Assume that $X$ has $r$ irreducible components $X_1,\ldots, X_r$. Then 
$H_1(M_X,\Z)=\Z^r$, where the generator $(0,\ldots,1,\ldots ,0) $ (1 on the 
place $i$, $1\leq i\leq r$) corresponds to an elementary oriented circle
``around $X_i$''. For each set of integers $\m :=(m_1,\ldots, m_r)\in \Z^r$
we define $\phi_{\m }:\Z^r\to \Z$ by $(x_1,\ldots, x_r)\mapsto
\sum_im_ix_i$. If $gcd_i(m_i)=\pm 1$ then
$$e_{X,\m }:\pi_1(M_X)\stackrel{ab}{\longrightarrow} H_1(M_X,\Z)\stackrel
{\phi_{\m }}{\longrightarrow}\Z \ \ \mbox{is onto}.$$

\vspace{2mm}

\noindent (II) \ Now, assume that $X=f^{-1}(0)$ for some polynomial $f$.
Then define $e_{X,f}$ by 
$e_{X,f}:=f_*: \pi _1(M_X) \to \pi _1 (\C^*)$.
In fact, this is a particular case of (I):
if $f=\prod_{i=1}^r g_i^{m_i}$ (where $g_i$ are irreducible with distinct
zero sets) then $e_{X,f}=e_{X,\m }$ for
$\m =(m_1,\ldots, m_r)$. 

By a similar argument as in the proof of \ref{t1}, if 
the generic fiber $F$ of $f$ is  connected then $e_{X,f}$ is onto. 
Therefore, in the sequel, in all our Alexander-module discussions associated 
with $f$, {\em we will assume that $F$ is connected}. 

\vspace{3mm}

Sometimes,  we will use the notations (I) resp. (II) to
remind the reader about the corresponding cases. 
In both cases (I) and (II), 
let $\F:=M_X^c$ be the $\Z$-cyclic covering 
associated with the kernel of $e_X$ (cf. also with \S 2).
Since $M_X$ has the homotopy type of a finite CW-complex, it follows that all the associated Alexander modules are of finite type over $\Lambda _A$
(but in general not over $A$).

Moreover, for a complex number $a \in \C^*$,
 we consider the rank one local system $L_a$ on $M_X$
defined by the composed map  $\pi _1(M_X)\stackrel{e_X}{
\longrightarrow} \Z\to\C^*$, where
the last map is defined by $1_\Z\mapsto a$. Obviously, if $a=1$, then $L_a=\C$.

Then, exactly as in  [Li3],  we have the following long exact sequence
\begin{equation*}
... \to H_k(\F,\C) \to H_k(\F,\C) \to H_k(M_X, L_a) \to  H_{k-1}(\F,\C) \to ...
\tag{2}
\end{equation*}
where the first morphism is multiplication by $t-a$.

\vspace{2mm}

In the next paragraph we summarize  the properties of 
$H_*(M_X,L_a)$ and the Alexander modules $H_*(\F,\C)$.

\subsection{Theorem}\label{3.6} {\it

(i)  $H_m(M_X,L_a)=0$ for any $a \in \C^*$ and $m > n+1$. 
  The Alexander modules $H_k(\F)$ are  trivial 
for $k>n+1$ and $H_{n+1}(\F,\C)$ is free (over $\Lambda$).

(ii) The $\Lambda$-rank of $H_{n+1}(\F,\C)$ equals $\dim H_{n+1}(M_X,L_a)$ 
for any generic $a $.

(iii) For any $q\leq n$, the $(t-a)$-torsion in $H_q(\F,\C)$ 
 can be non zero only when a is a root of unity.

(iv) Assume that all the Alexander modules $H_k(\F,\C)$
are torsion  for $k<n+1$. Denote by $N(a,k)$ the number of direct summands 
in the $(t-a)$-torsion part of $H_k(\F,\C)$. Then
$\dim H_k(M_X, L_a)=N(a,k)+N(a,k-1)$ for $k<n+1$ and
 $\dim H_{n+1}(M_X, L_a)=N(a,n) +|\chi (M_X)|$. 

\noindent Moreover, if either

(I) $X$ has only isolated singularities including at infinity, or

(II) $X$ is the fiber of a h-good polynomial,

\noindent then

(v)  $ \tilde{H}_k(\F)=0$ for $k<n$, and $H_n(\F,\C)$ is a
 $\Lambda$-torsion module.
Moreover, one has the isomorphisms of the Alexander $\Lambda_\Z$-modules
$\pi_n(M_X)=\pi_n(\F)=H_n(\F,\Z)$.}

\begin{proof}
For (i) and the first part of (ii) see \ref{fin}(2);
for (ii) use the exact sequence \ref{3.4}(2) and the fact that multiplication
 by  $t-a$ is injective, provided that $a$ is generic.
The vanishing of $(t-a)$-torsion for $a$ non root of unity follows from 
\ref{n7}(iv) below. (iv) is straightforward, if we notice that 
$\chi(M_X,L_a)=\chi(M_X)$ for any local system $L_a$.
Claim (v) in case (I) for $\m=(1,\ldots,1)$ is due to Libgober, see [Li2], 
Theorem 4.3; the case of arbitrary $\m$ follows similarly.
The case (II) is a consequence of \ref{con}.
\end{proof}

\subsection{Remarks.}\label{remarks} 

(1) Under the assumption (I) and (II) in \ref{3.6}, 
if $H_n(\F,\Z)$ is $\Z$-torsion free, then it has a finite rank over $\Z$ and the associated Alexander polynomial $\Delta (H_n(\F,\Z))$ is the characteristic polynomial $\Delta (h_n)$ as in \ref{3.3}. 
It follows
that $H_n(M_X,\Z)$ is finite iff $\Delta (h_n)(1) \ne 0$. Moreover, if 
$H_n(M_X,\Z)$ is finite then its order is exactly $|\Delta (h_n)(1)|$.

However, we do not know whether (I) or (II) in \ref{3.6}
 imply that $H_n(\F,\Z)$ is $\Z$-torsion free.

\vspace{1mm}

(2) \ The proof of \ref{3.6} is topological in both cases (I) and (II). 
It follows that the same proof  can be applied to coefficients in a finite field $K$, as soon as $|K|$ is large enough. Therefore, Theorem 
\ref{3.6} also
 holds in such a case (of course without any reference to roots of unity).

\vspace{1mm}

(3) \ Notice that \ref{3.6}(iv) is not true 
without the assumption about $\F$. In other words,
it is not true that  for a generic $a \in \C^*$  one has
$H_m(M_X,L_a)=0$ for $m \not= n+1$. Indeed, consider a polynomial function
$f:\C^{n+1} \to \C$ as a function $\C^{n+1+k} \to \C$ independent of the last $k$-variables and denote by $M_X'$ the corresponding complement. Then 
 $M_X' = M_X \times \C^k$. In particular, if $f$ is chosen such that $H_{n+1}(M_X,L_a) \neq 0$ for any generic $a$ (i.e. $\chi(M_X)\not=0$), we get counter-examples to the above claim.

\vspace{1mm}

(4) \ Additionally to the exact sequence \ref{3.4}(2) 
one has two more (both valid 
over $\Z$). First notice that the $\Z$-covering $\F\to M_X$ homotopically 
can be identified with the inclusion of the fiber
into the total space of a fibration $M_X\to K(\Z,1)=S^1$ (with fiber $\F$). 
Hence, the (homotopy)  Wang exact sequence of the covering $\F\to M_X$
gives:
$$... \to H_k(\F,\Z) \to H_k(\F,\Z)\to H_k(M_X,\Z) \to H_{k-1}(\F,\Z) \to ...$$
where the first morphism is multiplication by $t-1$.

Moreover, let $M_{X,e}$ ($e>0$) be the cyclic covering of $M_X$ of degree $e$
(i.e. the covering classified by the subgroup $e_X^{-1}(e\Z)$).
If $\Z$ denotes the transformation group of $\F$ above $M_X$, then 
 $\F$ is a $\Z$-covering of $M_{X,e}$ with transformation group $e\Z$. 
Hence, we have an exact sequence 
$$...\to H_k(\F,\Z)\to H_k(\F,\Z)\to H_k(M_{X,e},\Z)\to H_{k-1}(\F,\Z)\to ...$$
where the first morphism is multiplication by $t^e-1$. 
It follows that  there is a close relation between the homology of $M_{X,e}$ 
and the structure of the Alexander invariants $H_k(\F,\Z)$, see Example \ref{4.5} below.

\subsection{Example}\label{3.8} Assume that $f$ is a weighted homogeneous polynomial with respect to an integer set of weights (non necessarily strictly positive) such that $d=deg(f) \not=0$. Denote by $F$  the associated affine Milnor fiber. Then it is easy to see that $\F$ and $ F$ are homotopically equivalent,
and the covering transformation $h$ corresponds to the Milnor monodromy.
In particular $h_*^d=Id$. It follows that
$N(a,k)=\dim ker (h_*-a Id|H_k(\F,\C))$ is trivial for $a^d \not= 1$.
Moreover, in this case,
 due to the presence of $S^1$-actions, we have $\chi (M_X)=0$.

In the case of a central hyperplane arrangement D. Cohen [C] has shown that
for any rank one local system $L$ on $M_X$ one has
$\dim H_m(M_X,L) \leq \dim H_m(M_X,\C)$ for any $m$.

The example of the $A_1$-singularity $f:\C^3 \to \C$, $f=x^2+y^2+z^2$ , $m=2$ and $a=-1$,
 shows that this inequality is not true for a general homogeneous polynomial.

\section{Relations to Individual Monodromy Operators}

\subsection{}\label{n1} For the convenience of the reader we recall the
present set up. Let $f:\C^{n+1}\to \C$ be a polynomial   whose {\em generic 
fiber $F$ is connected}. Let $X=f^{-1}(0)$ be  a fixed  special fiber,
$M_X$ its complement and $e_X=f_*:\pi_1(M_X)\to\Z$ the epimorphism  whose 
kernel determines the $\Z$-cyclic covering $p:\F\to M_X$. Fix a generic fiber
$F\subset M_X$, and a base point $*\in \F$ with $p(*)=b_0\in F$. 
Then there is a 
unique lift (embedding) $s:F\to \F$  of $p$ over $F$  with $s(b_0)=*$. 
Let $h$ be the covering transformation of $\F$ corresponding to the generator
$1_\Z$. Sometimes we prefer to denote this Alexander $\Lambda_\Z$-module
 by $(H_*(\F),h_*)$. 

Now, fix an  orientation of $S^1$, and consider an arbitrary smooth map $
\lL:S^1\to \C\setminus B_f$
(we can even take  $\lL(1)=b_0$). 
Let $q:\lL^{-1}(f)\to S^1$ be the pull back of $f$ by $\lL$. 
Obviously, $q$ is a fiber
bundle over $S^1$ with fiber $F$. The map which covers $\lL$ is denoted
by $\bar{\lL}:\lL^{-1}(f)\to M_X$, hence $f\circ \bar{\lL}=\lL\circ q$.
We can consider the epimorphism $e_\lL:\pi_1(\lL^{-1}(f))\to \pi_1(S^1)=\Z$.
The associated $\Z$-covering has a total space isomorphic to $F\times \R$
and the covering transformation (corresponding to $1_\Z$)  can be identified
 with the geometric monodromy of $f$ associated with the oriented loop $\lL$.
This Alexander $\Lambda_\Z$-module will be denoted by 
$(H_*(F), T_{\lL,*})$. 

Next, we connect these two Alexander modules. First, assume that $\lL_*:
\Z\to\Z$ (i.e. $\lL_*:\pi_1(S^1)\to \pi_1(\C^*)$) satisfies $\lL_*(1)=1$. 
Then by \ref{t4}, $\pi_1(\bar{\lL}):\pi_1(\lL^{-1}(f))\to \pi_1(M_X)$ is 
onto, and $e_X\circ \pi_1(\bar{\lL})=e_\lL$. Therefore, by \ref{3.1} 
there exists a morphism of $\Lambda_\Z$-modules
$$\lL^c_{\Lambda,*}:(H_*(F), T_{\lL,*})\to (H_*(\F),h_*)$$
(in the sense that $h_*\circ \lL^c_{\Lambda,*}=\lL^c_{\Lambda,*}
\circ T_{\lL,*}$). 

More generally, assume that  $\lL_*:\Z\to \Z$ is multiplication with 
an integer $\ell$ 
(in other words, $\ell$ is the winding number of $\lL$ 
with respect to 0). Then 
we can replace $M_X$ by the $\Z/\ell\Z$-cyclic  covering $M_{X,|\ell|}$
(with projection $pr:M_{X,|\ell|}\to M_X$). By \ref{t4}, $\bar{\lL}:
\lL^{-1}(f)\to M_X$ can be lifted to  $\bar{\lL}':
\lL^{-1}(f)\to M_{X,|\ell|}$ (with $pr\circ \bar{\lL}'=\bar{\lL}$). Obviously,
$\F$ is the cyclic covering of $M_{X,|\ell|}$ with transformation group
$\ell\Z$, and $\pi_1(\bar{\lL}')$ is onto. Hence again we obtain a morphism 
of $\Lambda_\Z$-modules
\begin{equation*}
\lL^c_{\Lambda,*}:(H_*(F), T_{\lL,*})\to (H_*(\F),h_*^{\ell})
\tag{$\gamma$}
\end{equation*}
(i.e.  $h_*^{\ell}\circ \lL^c_{\Lambda,*}=\lL^c_{\Lambda,*}
\circ T_{\lL,*}$. Above, if
$\ell =0$ then $M_{X,|\ell|}$ is obviously $\F$ itself). 

The above constructions can be made compatible with some 
base points. (Since all the spaces are connected, this is not really
relevant, the details are left to the interested  reader). 

Notice that the target module  $(H_*(\F),h_*)$ is completely independent of
$\lL$. The main 
point is that even the map $H_*(F)\to H_*(\F)$, as a $\Z$-module,
is independent of $\lL$. Indeed, $\lL^{-1}(f)^c$ has the homotopy type of $F$
and $\bar{\lL}^c$ is up to homotopy the same as our fixed embedding 
$s:F\to \F$. The above discussion is summarized in the following
theorem.

\subsection{Theorem}\label{n2} {\em Assume that the generic fiber $F$ of $f$
is connected,  and $0\in B_f$. Fix an embedding $s:F\to \F$ as above. 
Let $\lL:S^1\to \C\setminus B_f$ be a smooth loop. Assume that 
$\lL_*:\pi_1(S^1)\to \pi_1(\C^*)$ is 
multiplication by $\ell=\ell(\lL)$. Then $s_*:H_*(F)\to H_*(\F)$ is compatible
with the Alexander $\Lambda_\Z$-module structures in the sense that
$s_*\circ T_{\lL,*}=h_*^\ell\circ s_*$. }

\vspace{2mm}

\noindent One has the following immediate consequence.

\subsection{Corollary}\label{n3} {\em Assume that two loops $\lL,\lL':S^1\to
\C\setminus B_f$ satisfy $\ell(\lL)=\ell(\lL')$ (i.e. they have the same 
winding number with respect to 0).
Then $T_{\lL,q}=T_{\lL',q}$ modulo $ker(s_q)$, for any $q\geq 0$.}

\vspace{2mm}

In other words, if the monodromy operators are ``very different'', 
then the image of  $s_*:H_*(F)\to H_*(\F)$  is forced to be ``small''.
Conversely, the non-vanishing of $im(s_*)$ imposes some compatibility 
restrictions on the monodromy operators. 
For a precise reinterpretation of \ref{n2} and \ref{n3}, see \ref{r4}.

The above theorem is optimal exactly when $s_q$ is onto in the non-trivial
cases $q\leq n$. As we will see,
this  is the case  e.g. for h-good polynomials (cf. \ref{n4} below). 
On the other hand, we cannot hope that for an  arbitrary polynomial $f$
the map 
$s_q:H_q(F)\to H_q(\F)$ ($q\leq n$)
is onto , since $H_q(F,\C)$ is a $\Lambda$-torsion
module, but $H_q(\F,\C)$ may have a non zero free part.
 Nevertheless, the next theorem shows that  $im(s_*)$ is 
as ``large as possible''.  Below  $\Delta_{\lL,*}(t)$ denotes 
the  characteristic polynomial of the monodromy operator $T_{\lL,*}$. 

\subsection{Theorem}\label{n7} {\em  Assume that $F$ 
is connected,  and $0\in B_f$, and $s:F\to \F$ is fixed as in \ref{n2}.
Let $({\mathcal T}_*,h_*)$ be the torsion part of the $\Lambda$-module 
$(H_*(\F,\C),h_*)$. Then:

(i) \ i$m(s_*)={\mathcal T}_*$.

(ii) Therefore, for any loop 
 $\lL:S^1\to \C\setminus B_f$ with  $\ell:=\ell(\lL)$ one has the following
epimorphism of $\Lambda$-modules:
$$s_*:(H_*(F,\C), T_{\lL,*})\to ({\mathcal T}_*,h_*^\ell).$$

(iii) In particular,  $\Delta((H_q(M_X^c,\C),h_q^\ell)(t)$  divides 
$\Delta_{\lL,q}(t) $ for any $q\geq 0$ (cf. with \ref{ellp}).

(iv)  Part (iii) applied for $T_{0,*}$ implies that the 
$(t-a)$-torsion of $H_*(M_X^c,\C)$ is zero if $a$ is not a root of unity.}

\begin{proof} Notice that part (i) and \ref{n1}($\gamma$) imply (ii), hence 
all the others as well. Notice that the statement of (i) is independent of the choice of any loop. In order to prove (i), we take a special loop, namely 
the oriented boundary of $D_0$. Then using the notations and the 
construction of \ref{con}, we deduce that 
$H_*(\F,F)=\oplus_{\bar{b}} H_*(\Z\times (T_{\bar{b}},F))$; in particular
this is a free $\Lambda$-module. Then the result follows from the 
homological exact sequence of the pair $(\F,F)$ (considered as a sequence of 
$\Lambda$-modules).
\end{proof}

\subsection{Remark}\label{ellp} If ${\mathcal M}=(H,h)$ is a 
torsion $\Lambda$-module
with Alexander polynomial $\Delta(H,h)(t)$ (i.e., if $h$ acts on $H$
with characteristic polynomial $\Delta(H,h)(t)$), then for any $\ell\geq 0$
the polynomial $\Delta(H,h)(t)$ determines 
$\Delta(H,h^\ell)(t)$ as follows.
For any polynomial  $P(t)=\prod_a(t-a)^{n_a}$ ($n_a\in \N$) write
$P(t)^{(\ell)}:=\prod_a(t-a^\ell)^{n_a}$. Then 
$\Delta(H,h^\ell)(t)=\Delta(H,h)(t)^{(\ell)}$.

\vspace{2mm}

Now, we will apply our general results for h-good polynomials. 
The next corollary follows directly from \ref{con} and \ref{3.2}.

\subsection{Corollary}\label{n4} {\em Assume that $f$ is a h-good polynomial,
$0\in B_f$, 
and take a smooth loop $\lL:S^1\to \C\setminus B_f$ with $\ell=\ell(\lL)$. 
Then there is an epimorphism of  $\Lambda_\Z$-modules
$$s_n:(H_n(F), T_{\lL,n})\to (\pi_n(M_X),h_n^{\ell}).$$
In particular,   $\Delta(\pi_n(M_X))(t)^{(\ell)}$ 
divides $\Delta_{\lL,n}(t) $.}

\subsection{The main examples. 
The monodromy around the origin and at infinity}\label{n5} 

To any polynomial $f$ (and to a distinguished atypical value $0\in B_f$)
 one can associate two distinguished monodromy operators:

(i) the local monodromy of $f$ at 0, namely $T_{0,*}: H_*(F) \to H_*(F)$ 
provided by a loop $\lL$ going  around the bifurcation point 0 on 
the boundary of a small disc $D_0$ containing no other bifurcation 
points inside (with positive orientation);

(ii) the monodromy at infinity of $f$, namely $T_{\infty,*}: H_*(F) \to 
H_*(F)$ provided by a loop $\lL$ going around the boundary of a large disc 
$D_{\infty}$ containing all the bifurcation points $B_f$ inside.

They provide two (Alexander) $\Lambda$-module structures on 
$H_*(F,\C)$ denoted by  $A_*(f,0)$), respectively  $A_*(f,\infty)$). 
The corresponding Alexander polynomials, or equivalently, the characteristic
 polynomials of the monodromy operators $T_{0,*}$ and $T_{\infty,*}$,
are denoted by $\Delta(T_{0,*})$ resp. $\Delta(T_{\infty,*})$. 
Note that the Alexander module $A_*(f,0)$ (resp. $A_*(f,\infty)$) encodes
exactly the Jordan structure of $T_{0,*}$ 
(resp. $T_{\infty.*}$) and that this Jordan structure was studied 
in several papers, e.g. [ACD], [D2], [DN1-2], [DS], [Do], [GN1-3], [Sa].

Obviously, in both cases $\ell=1$. Therefore, \ref{n7} guarantees that
for any $f$ with $F$ connected, and for any $q\geq 0$,
the Alexander polynomial 
\begin{equation*}
\mbox{
$\Delta ( H_q(M_X^c, \C))$ divides the 
characteristic polynomials $\Delta (T_{0,q})$ and $\Delta (T_{\infty,q})$.}
\tag{1}
\end{equation*}
This is a generalization of [Li2] (4.3) to arbitrary polynomials and $q$. 
Let us explain more precisely
the relation between this  divisibility result (1)
 and the divisibility results in [Li2]. 

First consider the local monodromy  $T_{0,*}$.
In [Li2],  $X$ has only isolated singularities including at infinity.
In that case only the case $q=n$ is relevant (cf. \ref{3.6}). 
It is known
that by the localization of the monodromy (see e.g. [NN]) one has
\begin{equation*}
 \Delta (T_{0,n})(t)=(t-1)^k\cdot \prod \Delta _i(t),
\tag{2}
\end{equation*}
where $\Delta _i$ are the characteristic polynomials of the local 
monodromies associated to the isolated singularities on $V$, 
and $k =b_n(F) - \sum deg(\Delta _i)$. In fact,
 the singularities at infinity (i.e. on $V \cap H$)
 should be treated slightly different, 
as explained in [Li2], [LiT]. 
Moreover, in general $k \geq 0$,  and $k=0$ iff $B_f=\{0\}$.
Note that (1) (for $q=n$)  and (2) give a similar result to Theorem (4.3)
 in [Li2], yielding in addition a precise value for $k$.

The discussion for the monodromy at infinity $T_{\infty,*}$ is more involved.
It was shown by Neumann and Norbury [NN] that the total space of the fibration
($*$)\ 
$ f: f^{-1}(S^1_r) \to S^1_r$ for $r\gg 0$ (which provides $T_{\infty,*}$)
can be embedded in a natural way as an open subset of $S^{2n+1}\setminus
f^{-1}(0)$, where $S^{2n+1} $ is a large sphere in $\C^{n+1}$.
Moreover, it was shown in [NZ] that for an M-tame polynomial this fibration
($*$) is equivalent to the Milnor fibration at infinity
$\phi : S^{2n+1} \setminus f^{-1}(0) \to S^1$, $\phi (x)=f(x)/|f(x)|$.

Note that in general $S^{2n+1} \setminus f^{-1}(0)$ is not the total space of a fibration over the circle, or, even when it is, it may happen that the corresponding fiber is not $F$,
 see the case of semi-tame polynomials in P\u aunescu-Zaharia [PZ].

In [Li2],  Libgober considers an infinite cyclic covering $U_{\infty}$ of the knot complement $S^{2n+1} \setminus f^{-1}(0)$ and takes the associated Alexander module $H_n(U_{\infty},\C)$ as the Alexander module at infinity for $f$.
From our discussion above it seems that, 
in general,  one cannot hope to identify easily the structure of this module 
$H_n(U_{\infty},\C)$. However, in the case of $M$-tame polynomials,
the module $H_n(U_{\infty},\C)$ is exactly our $A_n(f,\infty)$.

\subsection{Example}\label{4.5} Recall the situation described in \ref{defect}.
Namely, let $f$ be a polynomial such that there exists a system of weights $w$ with the top degree form $f_e$ defining an isolated singularity at the origin.
The monodromy at infinity of such a polynomial coincides to the monodromy of the singularity $f_e=0$, 
in particular $T_{\infty,n}$ is semisimple and all the eigenvalues are 
$e$-roots of unity.
It follows that in the second exact sequence in \ref{remarks}(4), 
for  $k=n$, the first morphism is trivial. Hence (for $n>1$) 
$$\pi _n(M_X) \otimes \C =   H_n(\F,\C) =H_n(M_{X,e},\C)=H_n(F',\C)$$
are isomorphic $\Lambda$-modules. This results should be compared to Corollary (4.9) in [Li2].

As a concrete example, let $f:\C^4 \to \C$ be the tame polynomial considered in Example \ref{2.10}. The above discussion and \ref{n5}(1) gives:
$$\pi _3(M_X) \otimes \C =   H_3(\F,\C) =H_3(M_{X,3},\C)=H_3(F',\C)= (\Lambda /(t-1))^5.$$

\subsection{Remarks}\label{4.6}
If one wants to determine the homology groups of $M_X$, either one needs
some information about the  ``tubes''
$T_{\bar{b}}$ for $\bar{b}\in B_f\setminus \{0\}$, or (using \ref{2.7}) 
one needs to know the behaviour at infinity of $X$; both rather subtle 
problems. Therefore, it is rather
surprising that, in some cases,  all these information is 
carried by only one monodromy operator $T_{0,n}$.

Here we present the case of h-good polynomials:
we describe completely  $H_*(M_X,\Z)$ in terms of $T_{0,n}$.

Let ${\mathcal V_0} \subset H_n(F,\Z)$ be the subgroup of vanishing cycles 
at 0 corresponding to a choice of a star in $\C$ as in [DN1]. 
It follows as in [{\em loc.\,cit.}] that the morphism $T_{0,n}-Id$ induces a 
``variation'' morphism $V: H_n(F,\Z) \to {\mathcal V}_0$ and, 
by restriction to ${\mathcal V}_0$, a morphism $V_0: {\mathcal V}_0 
\to {\mathcal V}_0$. Using the definition of h-good polynomials, 
the connectivity \ref{p23} of $F$, and the Wang sequence associated to 
$T_0^*$ (cf. also with \ref{con2}), one can prove the following.

\vspace{2mm}

 The homology groups of $M_X$ are trivial except possibly for:

\noindent {\bf  Case $n>1$.} \

(i) $H_0(M_X,\Z)=H_1(M_X,\Z)=\Z$,

(ii) $H_{n+1}(M_X,\Z)$ is $\Z$-torsion free of rank $b_n(F)-\rank(V)$ and

(iii) $ H_{n}(M_X,\Z)=\coker(V)$.

\noindent {\bf  Case $n=1$.}\

(i') $H_0(M_X,\Z)=\Z$,

(ii') $H_{2}(M_X,\Z)$ is $\Z$-torsion free of rank $b_n(F)-\rank(V)$,

(iii') $ H_{1}(M_X,\Z)=\coker(V) +\Z$, a $\Z$-torsion free group of rank 
equal to the number of irreducible components of $X$.

\vspace{2mm}

Note that for $n>1$, $ H_{n}(M_X,\Z)$ is a finite group if $V_0$ is injective.
This happens exactly when, in the notation from \ref{n5}(2), 
one has $\prod \Delta _i(1) \ne 0$. Moreover, the epimorphism 
$\coker(V_0) \to \coker(V)$ implies that the order 
$|H_{n}(M_X,\Z)|$ divides $|\prod \Delta _i(1)|=|\coker(V_0)|$.

Let $g:\C^{n+1} \times \C \to \C$ be the $d$-suspension of the polynomial $f$, namely $g(x,y)=f(x)-y^d$. Let $Y=g^{-1}(0)$. Then writing the Gysin sequences in homology associated to a smooth divisor $D$ in a complex manifold $Z$ for the pairs $(Z,D)=(M_Y,M_X), \
(\C^{n+1} \times \C^*, M_{X,d})$ and resp. $(M_X \times \C^*, graph(f))$,
 and comparing the associated morphisms, we get for all $q>0$
the  exact sequence
\begin{equation*}
 0 \to H_{q}(M_X,\Z) \to  H_{q}(M_{X,d},\Z)\to  H_{q+1}(M_Y,\Z) \to 0.
\tag{3.q}
\end{equation*}
The exact sequence $(3.n+1)$ is  split since the last group in it is free 
according to (ii) above. The exact sequence $(3.n)$ is not split, as can be 
seen in the case $n=1,f=x_1x_2$, $d=3$ when we get 
$0 \to \Z^2 \to \Z^2 \to \Z/3\Z \to 0$. 
This example shows the difficulty of the question \ref{2.12}.

Finally, assume that $n>1$ and $\prod_i \Delta _i(\alpha ^k) \ne 0$
for $\alpha =exp(2\pi i/d)$ and for any $k \in \Z$ (cf. \ref{n5}(2)).
Then one has:

\vspace{2mm}

(a) all the groups $ H_n(M_X,\Z),\ H_n(M_{X,d},\Z)$ and $ H_{n+1}(M_Y,\Z)$ are finite; and

(b) the order $|H_n(M_{X,d},\Z)|$ divides the product $|\prod _{1\leq k\leq d}
 \prod_i \Delta _i(\alpha ^k)|$.

\vspace{2mm}

The proof of these claims follows from the exact sequence (3.q) 
 and the property (iii) above once we know how to compute the variation associated to the special fiber $Y$. This, in turn,
 is explained in [DN2]. 
Note that the claim (b) is similar to Theorem 3 in [Li0].

\section{ Relations to Monodromy Representation}

\subsection{}\label{r1} The results of the previous section already suggest 
(see e.g. \ref{n3})
that one can obtain finer results about the Alexander modules if one takes 
the whole monodromy representation instead of individual monodromy operators. 
 The main message of this section is that from
the monodromy representation of $f$  one can construct a universal Alexander
module which, in some sense, dominates all the Alexander modules 
associated with (all) the fibers  of $f$. 

Since the case of h-good polynomials with  all the involved numerical
invariants (cf. \ref{5.1} and \ref{r2})
represents a special interest, we start our detailed
discussion with this case. But, thanks to the general results of the 
previous sections, the next  constructions and factorization 
phenomenon described in the h-good polynomial 
case, can be repeated word by word in the general case. The general
result will be formulated at the end of the section in \ref{gen}. 

We start with a h-good polynomial. 
With the notation of \S 2, let $S=\C \setminus B_f$,  $E=f^{-1}(S)$
and $g=|B_f|$. Then the locally trivial fibration $f:E \to S$ induces a monodromy representation $\rho: G \to Aut (\L)$, where $G=\pi _1(S,b_0)$ is a free 
group on $g$ generators,
$b_0\in S$ is a base point,  and $\L=H_n(F,\Z)$ with $F=f^{-1}(b_0)$. 
For each $b\in B_f$ write $F_b=f^{-1}(b)$. 
Let $\gamma _i$ denote an elementary loop around $b_i \in B_f$ and 
$m_i=\rho ( \gamma _i)$ be the corresponding monodromy operators. With a natural choice for $\{\gamma _i\}_i$ one has $m_1 \cdot m_2 \cdot ... \cdot 
m_g =T_{\infty,n}$, see [DN1].

For any $H$-module $\M$ of a group $H$, we 
 denote by $\M_H$  the group of coinvariants, namely the quotient of $\M$ by 
the subgroup spanned by all elements $h \cdot m -m$ for $h \in H$ and 
$m \in \M$, see Brown [Br].
We denote by $b_k^c(Y)$ the $\C$-dimension of the $k^{th}$-cohomology space $H^k_c(Y,\C)$ of $Y$ with compact supports. 

We start by recalling how the $G$-module $\L$ determines the homology of the 
space $E$. 

\subsection{Proposition}\label{5.1} {\it The reduced homology groups 
${\tilde H}_k(E,\Z)$ are trivial except at most for $k=1$, $k=n$ and $k=n+1$. For these values of $k$ one has the following.

(i) For $n=1$ one has $H_2(E,\Z)=H_1(G,\L)$ and an exact sequence of groups
$$ 0 \to \L_G \to H_1(E,\Z) \to \Z^g \to 0.$$
In particular, $\L_G$ is a free $\Z$-module with} $\rank (\L_G) = 
\sum _{b \in B_f}(n(F_b)-1)$, {\em
where $n(Y)$ denotes the number of irreducible components of a curve $Y$.

(ii) For $n>1$ one has $H_1(E,\Z)=\Z^g$, $H_n(E,\Z)=H_0(G,\L)=\L_G$ and 
$H_{n+1}(E,\Z)=H_1(G,\L)$. In particular, } $\rank (\L_G)= 
\sum _{b \in B_f} b_{n+1}^c(F_b)=
\sum _{b \in B_f} b_{n+1}(T_b,\partial T_b)$.

\vspace{2mm}

\noindent {\em Proof.} The result follows from the Leray spectral sequence in homology of the fibration $F \to E \to S$ and basic facts on group homology, see Brown [Br]. The claim about the rank of $\L_G$ follows from the long exact sequence
$$\hspace{2.5cm}...  \to H^k_c(E) \to H^k_c(\C^{n+1}) \to H^k_c(\cup F_b) 
\to H^{k+1}_c(E) \to ...$$

\vspace{1mm}

Let $H=[G,G]=G'$ be the commutator of $G$ and $S' \to S$ be the corresponding covering space. Let $f':E' \to S'$ be the fibration (with fiber $F$) 
obtained from the fibration $f:E \to S$ by pull-back. Then the monodromy of the fibration $f'$ corresponds exactly to the $H$-module $\L$ obtained by restriction of $\rho$ to $H$. On the other hand, we can regard $E' \to E$ as being the covering space
corresponding to the kernel of the composition $\pi _1(E) \to \pi _1(S)=G \to G/H=\Z^g \to 0$. It follows that the deck transformation group of $E' \to E$ is $\Z^g$ and hence we can regard $H_n(E',R)$ as a $\Lambda _{R,g}$-module, where 
$\Lambda _{R,g} =R[\Z^g]$ is a Laurent polynomial ring in $g$ indeterminates $t_1,...,t_g$.
As before, when $R=\C$ we simply write $\Lambda _{g}$.

To state the result similar to \ref{r1} for the fibration $f':E' \to S'$, note that $S'=\R$ for $g=1$ and $S'$ is homotopy equivalent to a bouquet of infinitely many $S^1$'s for $g>1$.

\subsection{Proposition}\label{r2} {\it The reduced homology groups ${\tilde H}_k(E',\Z)$ are trivial except at most for $k=1$, $k=n$ and $k=n+1$. For these values of $k$ one has the following.

(i) For $n=1$ one has $H_2(E',\Z)=H_1(H,\L)$ and an exact sequence of groups
$$ 0 \to \L_H \to H_1(E',\Z) \to H_1(H,\Z) \to 0.$$

(ii) For $n>1$ one has $H_1(E',\Z)=0$ for $g=1$ and $H_1(E',\Z)=H_1(H,\Z)$ for $g>1$, $H_n(E',\Z)=H_0(H,\L)=\L_H$ and $H_{n+1}(E',\Z)=H_1(H,\L)$.}

\vspace{2mm}

\noindent 
Using the description of the $K(H,1)$ and of the associated chain complex given in  [Li4], (1.2.2.1),  it follows that
$H_1(H,\Z)=G'/G''$ is a submodule of $\Lambda _{\Z,g}^g$ and hence 
$H_1(H,\Z)$ is $\Lambda_{\Z,g}$-torsion free. 
This implies that in both cases (i) and (ii) in 
\ref{r2},  we have $\L_H =\mbox{Tors}(H_n(E',\Z))$ (as a 
$\Lambda_{\Z,g}$-module).

In the sequel we denote the  $\Lambda_{\Z,g}$-module  $\L_H$ 
 by $M(f)$, and we call it the
{\it global Alexander module of the polynomial $f$}.

\subsection{Remark}\label{r3} 
The global Alexander module of the polynomial $f$ can be regarded as a 
commutative version of the monodromy representation $\rho$.
Notice also that using Brown [Br], Exercise 3, p.35, it follows that
$M(f)_{\Z^g}=\L_G$.

\vspace{2mm}

Assume now that $0 \in B_f$. We will construct a new $\Lambda_\Z$-module 
$M(f,0)$ out of the monodromy representation. The inclusion $S \to \C^*$ 
at $\pi_1$-level gives rise to a projection $p_0: G \to \Z$. 
Let $K_0$ denote the kernel of this projection. Then $H \subset K_0 $ and hence we have a tower of covering spaces $S' \to S^0 \to S$. Let $E^0 \to S^0$ be the fibration induced from $E \to S$ by pull-back.
In this way we get a second tower of covering spaces, namely 
$E' \to E^0 \to E$.

We have, exactly as in the proof of \ref{r1}, 
the following isomorphisms of $\Lambda$-modules:
$$\mbox{Tors}(H_n(E^0,\Z))=\L_{K_0}=(\L_H)_{K_0/H}.$$
Here $K_0/H =\Z^{g-1}$ with generators corresponding to elementary loops in $\C$ around points in $B_f$ different from 0. Moreover the $\Z=G/K_0$-action on
$\mbox{Tors}(H_n(E^0,\Z))$ is induced by the monodromy operator $T_{0,n}$. 

We denote this 
$\Lambda_\Z$-module $\L_{K_0}$ by $M(f,0)$,
 and we call it  the {\it local Alexander module of $f$ at 0.}
(Clearly,  similar local Alexander module can be defined for any $b\in B_f$.)

The above isomorphisms show that the local  Alexander module $M(f,0)$ of $f$ at 0 can be computed from the global  Alexander module $M(f)$ of $f$. In this sense, the module $M(f)$ is universal, i.e. contains all the information about
the local Alexander modules associated to all the special fibers of $f$.

%Both modules $M(f)$ and $M(f,0)$ can be considered with coefficients in a commutative ring $R$ and in this case they will be denoted by $M(f;R)$ and $M(f,0;R)$. The coefficients $R$ are omitted when $R=\Z$ or when the choice for $R$ is clear from the context.

The usefulness of this new Alexander module comes from the fact that it can be calculated using the monodromy representation and gives another approximation for the Alexander module 
$\pi_n(M_X)=H_n(\F,\Z)$. Before we formulate this statement, let us
reinterpret \ref{n2} and \ref{n3}.

\subsection{Theorem \ref{n2} revisited}\label{r4} Let us explain the meaning 
of \ref{n2} in the language of the present section. Clearly, $\L=
H_n(F,\Z)$ is a 
$G$-module, and $H_n(\F,\Z)$ has a cyclic action generated by $h_n$. The 
map $p_0:G\to \Z$ can be identified with $[\lL]
\mapsto \ell(\lL)$ considered in 
\ref{n2}. Therefore, if we consider both $\L$ and $H_n(\F,\Z)$ as $G$-modules
(the last one via $p_0$), then \ref{n2} says that $s_n$ is a morphism of 
$G$-modules. 

In other words, the complicated monodromy representation (i.e. the $G$-module
$\L$), when it is mapped via $s_n$ into $H_n(\F,\Z)$, it is collapsed
into  a modest cyclic action. Since the action in the target is abelian, 
this already 
shows that $s_n:\L\to H_n(\F,\Z)$ has a factorization through $\L_H=M(f).$

Notice that $K_0=\ker(p_0)$ constitutes of loops (with base points)
$\lL$  with $\ell(\lL)=0$. Corollary \ref{n3} applied for such a loop $\lL$
 and for the trivial loop guarantees that $\rho(\lL)m-m\in \ker(s_n)$
for any $m$. In particular, 
$s_n:\L\to H_n(\F,\Z)$ has the following factorization of $G$-modules:
\begin{equation*}
\L\to\L_H\to \L_{K_0}\to H_n(\F,\Z).
\tag{1}
\end{equation*}

\subsection{Corollary}\label{r6} {\em 
Assume that $f:\C^{n+1} \to \C$ is a h-good polynomial. 
Then $s_n:H_n(F,\Z)\to H_n(\F,\Z)$ induces an
epimorphism $M(f,0) \to  \pi_n(M_X)$ of $\Lambda _{\Z}$-modules.}
\begin{proof}
Since $s_n$ is epimorphism (cf. \ref{n4}), the result follows from (1) above.
\end{proof}
 
\subsection{Remark}\label{r8} Notice that any $\lL$  with $\ell(\lL)=+1$
induces the same operator $\overline{\rho}([\lL])$ acting on $\L_{K_0}$;
and this operator is the positive generator of the cyclic action on
$\L_{K_0}$. 
E.g., one can take $\overline{T}_{0,n}$, or $\overline{T}_{\infty,n}$ as well,
depending which one is easier to compute. 
We write  $M(f,0)=(\L_{K_0},\overline{T}_{0,n})$.
Then, for an arbitrary 
 $[\lL]\in G$ with  $\ell:=\ell(\lL)=p_0([\lL])$ one has the 
$\Lambda_\Z $-module epimorphisms:
$$(H_n(F,\Z),\rho(\lL)) \to (\L_{K_0},\overline{T}_0^\ell) \to  (\pi_n(M_X),h_n^\ell).
$$
Evidently, this provides the divisibilities of the corresponding Alexander
(or characteristic) polynomials. 

\subsection{Example}\label{r7} Assume we are in the situation of Example 
\ref{2.3} with $n>0$ even. Then it is easy to see that
$\L_G=\Z/2\Z$ and
$\L=\L_H=M(f)=M(f,0)=\pi_n(M_X)=\Lambda _{\Z}/(t-1)$.

\vspace{2mm}

Using the general result \ref{n2} and \ref{n7}, one can verify easily that
the above factorization (1)  is valid for arbitrary polynomials as well. 

\subsection{Theorem}\label{gen} {\em Let $f$ be an arbitrary polynomial 
with $F$ connected. For any $q\geq 0$, consider $\L_q:=H_q(F,\Z)$ as 
a $G=\pi_1(S,b_0)$-module provided by the  monodromy representation.
Define the global Alexander $\Lambda_{\Z,g}$-module by $(\L_q)_H$, 
and the local Alexander module associated with the bifurcation point $0\in B_f$
by $(\L_q)_{K_0}$. If we consider $H_q(\F,\Z)$ as a $G$-module via $p_0$, then 
$s_q:\L_q\to H_q(\F,\Z)$ has the following factorization  of $G$-modules:
\begin{equation*}
s_q: \L_q\to(\L_q)_H\to (\L_q)_{K_0}\to H_q(\F,\Z).
\end{equation*}
If one tensor this tower by $\C$, then the last term $H_q(\F,\C)$ can be
replaces by ${\mathcal T}_q$, being the image of $s_q$. }

\section{An Example}

Let $f:\C^2 \to \C$ be the polynomial $f=x+x^2y^2+x^2y^3$. Then $B_f=\{b_1,b_2\}$ with $b_1=-27/16$ and $b_2=0$. The fiber $F_{b_1}$ is irreducible, has a node as a 
singularity and is regular at infinity. 
On the other hand, the fiber $F_0=F_{b_2}$ is smooth, has two irreducible 
components, one a copy of $\C$ the other $\C \setminus \{0,-1\}$,
 and has a singularity at infinity with a Milnor number equal to 3. 
It follows that $b_1(F)=4$ and the Jordan normal form for the monodromy operators $m_1, m_2$ and $ T_{\infty}$ was obtained in [BM]:

$$m_1 \approx \begin{pmatrix}
 1&1&0&0 \\ 0& 1&0&0 \\ 0&0& 1 &0 \\ 0& 0& 0& 1 \end{pmatrix}
 \ \ m_2 \approx   \begin{pmatrix} 1&0&0&0 \\ 0&1&0&0 \\ 0&0&j&0 \\ 0&0&0&j^2 
\end{pmatrix} \ \
T_{\infty} \approx \begin{pmatrix} -1&1&0&0 \\ 0& -1&0&0 \\ 0&0& 1 &0 \\
 0& 0& 0& 1 \end{pmatrix}$$
with $j=exp(2 \pi i/3)$.

Let $g:\C^2 \to \C$ be given by $g=u+u^2v$. Then $B_g=\{0\}$, the generic fiber is homotopy equivalent to $S^1$ and the corresponding monodromy $m$ is the identity. Consider the polynomial $h:\C^4 \to \C$ given by $h(x,y,u,v)=f(x,y)+g(u,v).$

It follows from [DN2] that $B_h \subset B_f$.
Notice that $f$ and $g$ are not good, but  both are h-good. 
 The fact that $h$ is also a h-good polynomial follows from [DN2], 
Corollary (4.4), which basically says that the 
``Thom-Sebastiani sum'' of two h-good polynomials is h-good. 

Our goal is to determine 
the various  Alexander modules associated with $h$ and to its special fibers. 
Since the information we have on the monodromy representation of $f$ is over 
$\C$,  we choose this coefficient ring.

The generic fiber of $h$ is the join $(\vee_4 S^1)*S^1$, hence it is 
$\vee_4S^3$. Since $m=Id$, [DN2] guarantees  that
the monodromy representation of $h$ can be identified to that 
of  $f$. Using the above Jordan forms:
\begin{equation*}
A_3(h,b_1)=\Lambda/(t-1) \oplus \Lambda/(t-1)\oplus \Lambda/(t-1)^2.\tag{1}
\end{equation*}
\begin{equation*}
A_3(h,0)=\Lambda/(t-1) \oplus \Lambda/(t-1)\oplus \Lambda/(t-j) \oplus \Lambda/(t-j^2);
\tag{2}
\end{equation*}
\begin{equation*}
A_3(h,\infty)=\Lambda/(t-1) \oplus \Lambda/(t-1)\oplus \Lambda/(t+1)^2;
\tag{3}
\end{equation*}
Now, we consider the space $M_X$ corresponding to the two special fibers.
First, let $X=h^{-1}(0)$. Then $X$ is smooth and applying Theorem (4.7), (ii) 
in [DN2],  we get that $X$ has the homotopy type of $S^2 \vee S^3 \vee S^3.$
It follows that $b_3(M_X)=b_2(X)=1$. 
Using  \ref{n5}(1)  and \ref{3.6}(iv) one has $N(1,3)=1$, hence:
\begin{equation*}
\pi _3(M_X) \otimes \C =   H_3(\F,\C) =H_3(M^c_{X},\C)= \Lambda /(t-1).
\tag{4}
\end{equation*}
Next, let $X=h^{-1}(b_1)$. Note that $X$ is again smooth but we can no longer apply Theorem (4.7) in [DN2] since $F_{b_1}$ is not smooth.
Using the equality $ \chi (X) = \chi _c (X)$ we get
$b_3^c(X) -b_4^c(X)=3$. Moreover, it is known that $b_3^c(X)= 
\dim\coker (m_1-1)=3$, see [ACD] or [DN1]. It follows that $b_4^c(X)=0$. The exact sequence 
$$...  \to H^k_c(M_X) \to H^k_c(\C^{n+1}) \to H^k_c(X) \to H^{k+1}_c(M_X) \to ...$$
 and duality for $M_X$ imply that $b_3(M_X)=b_4^c(X)=0$. Moreover,
from (1) follows that only $(t-1)$-torsion is possible. Hence,
via \ref{3.6}(iv),   $b_3(M_X)=0$ implies that
\begin{equation*}
\pi _3(M_X) \otimes \C =   H_3(\F,\C) =H_3(M^c_{X},\C)= 0.
\tag{5}
\end{equation*}
Our next aim is to compute the global Alexander module $M(h)$. This is done by
using the partial information we have on the monodromy representation 
$\rho: G \to Aut (\L)$, where $G$ is a free group on two generators, and $\L$
 is the third homology of the generic fiber of $h$ with $\C$ coefficients
(i.e.,  to simplify notation, we denote by $\L$ the complexification 
$\L \otimes _{\Z} \C$).
In terms of a special basis $e_1,e_2,e_3,e_4$ for $\L$ as in [DN1], (2.5), 
we can write the monodromy operators $m_1$ and $m_2$  of $h$ in the form
$$m_1 = \begin{pmatrix} 1&a&b&c \\ 0& 1&0&0 \\ 0&0& 1 &0 \\ 0& 0& 0& 1 
\end{pmatrix}
 \ \ m_2= \begin{pmatrix} 1&0&0&0 \\ \alpha &1&0&0 \\ \beta &0&j&0 \\
 \gamma &0&0&j^2 \end{pmatrix}. $$
Checking that $m_1m_2$ is conjugate to the Jordan normal form for $m_{\infty}$ given above implies $a \alpha =0$ and $bc \not= 0$. By an obvious change of base we may take $b=c=1$ and then $\beta$ and $\gamma$ are determined by the equations $\beta +\gamma =-1, \beta j + \gamma j^2=2$.

Let $C=[m_1,m_2]=m_1m_2m_1^{-1}m_2^{-1}$. Then $C \in H$ and a direct computation shows that $v_1=(C-Id)(e_1)=-(\alpha e_2 +\beta e_3 +\gamma e_4)$.
Let $\L_0$ be the vector subspace in $\L$ spanned by all the vectors $h(v)-v$ for $h \in H$ and $v \in \L$. It follows that

(i) $v_1 \in \L_0$ and

(ii) $\L_0$ is a $G$-invariant subspace of $\L$ (this property being always true).

It follows that we have to discuss two cases.

Case 1. ($\alpha \ne 0$) Then the vectors $v_1$, $m_2v_1$ and $m_2^2v_1$ span the same subspace in $\L$ as the vectors $e_2,e_3,e_4$. Moreover $m_1e_3=e_1+e_3 \in \L_0$. Therefore $\L=\L_0$ and hence $\L_H=\L/\L_0=0$. 
But this is a contradiction since we have epimorphisms $M(h)=\L_H \to M(h,0) \to H_3(M^c_{X},\C)= \Lambda /(t-1)$ by \ref{r6} and (4)  above.

Case 2. ($\alpha = 0$) As above one shows that $\L_0$ is spanned by $e_1,e_3, e_4$ and hence $\L_H=\C$ with a trivial $\Z^2$-action. This implies the following.

\subsection{Proposition}\label{6.3} {\it For the polynomial $h:\C^4 \to \C$ described above one has the following Alexander modules.

(i) For the fiber $X=h^{-1}(0)$ one has $\pi _3(M_X) \otimes \C =H_3(M^c_{X},\C)= M(f,0)=\Lambda /(t-1).$

(ii) For the fiber $Y=h^{-1}(b_1)$ one has $\pi _3(M_Y) \otimes \C =H_3(M^c_{Y},\C)= 0.$ Moreover in this case $M(f,b_1)=\Lambda /(t-1) \ne H_3(M^c_{Y},\C)= 0$.

(iii) The global Alexander module $M(h)$ is isomorphic to $\Lambda _2/(t_1-1,t_2-1)$.}

\vspace{2mm} 

\noindent Note that:

(i) we succeeded to determine the above data in spite of the fact that
the monodromy representation of $h$ (over $\C$) is not completely determined 
(the value of $a$ is unknown);

(ii) in the case of $Y$, we have $M(h,b_1) \ne H_3(M_Y^c,\C)$,
in particular we cannot expect isomorphism in \ref{r6}.
Nevertheless, the approximation of $H_3(M_Y^c,\C)=0$ given by $M(h,b_1)$ is better than that given by $A_3(h,b_1)$ since $ \dim_{\C}M(h,b_1)=1$ while $ \dim_{\C}A_3(h,b_1)=4.$

\bigskip

Laboratoire de Math\'ematiques Pures de Bordeaux

Universit\'e Bordeaux I

33405 Talence Cedex, FRANCE

email: dimca@math.u-bordeaux.fr
\bigskip

Department of Mathematics,

Ohio State  University,

Columbus, Ohio 43210, USA

email: nemethi@math.ohio-state.edu

\end{document}